\documentclass[12pt]{amsart}
\usepackage{amsmath,latexsym,amscd,amsbsy,amssymb,amsfonts,amsthm,fleqn,leqno,
euscript, graphicx, texdraw}
\usepackage{pb-diagram}
\numberwithin{equation}{section}

\newtheorem{thm}{Theorem}[section]

\newtheorem{lem}[thm]{Lemma}
\newtheorem{cor}[thm]{Corollary}

\theoremstyle{definition}

\theoremstyle{remark}

\begin{document}

\title[Probability measures and Milyutin maps between metric spaces]
{Probability measures and Milyutin maps between metric spaces}

\author{Vesko Valov}
\address{Department of Computer Science and Mathematics,
Nipissing University, 100 College Drive, P.O. Box 5002, North Bay,
ON, P1B 8L7, Canada} \email{veskov@nipissingu.ca}

\date{}
\thanks{The author was partially supported by NSERC
Grant 261914-08.}

\keywords{Probability measures, support of measures, Milyutin maps,
metric spaces} 

\begin{abstract}
We prove that the functor $\Hat{P}$ of Radon probability measures
transforms any open map between completely metrizable spaces into a
soft map. This result is applied to establish some properties of
Milyutin maps between completely metrizable spaces.
\end{abstract}
\maketitle\markboth{}{}



\section{Introduction}

In this paper we deal with metrizable spaces and continuous maps. By
a (complete) space we mean a (completely) metrizable space, and by a
measure a probability Radon measure. Recall that a measure $\mu$ on
$X$ is said to be:

\begin{itemize}
\item {\em probability} if $\mu(X)=1$;
\item {\em Radon} if $\mu(B)=\sup\{\mu(K):K\subset B\hbox{~}\mbox{and $K$
is compact}\}$ for any Borel set $B\subset X$;
\end{itemize}

The support $supp\hbox{~}\mu$ of a measure $\mu$ is the intersection
of all closed subsets $A$ of $X$ with $\mu(A)=\mu(X)$. It is well
known that the support of any measure is non-empty and separable.

Everywhere below $\Hat{P}(X)$ stands for the space of all
probability Radon measures on $X$ equipped with the weak topology
with respect to $C^*(X)$. Here, $C^*(X)$ is the space of bounded
continuous functions on $X$ with the uniform convergence topology.
According to \cite{b1}, $\Hat{P}$ is a functor in the category of
metrizable spaces and continuous maps. In particular, for any map
$f\colon X\to Y$ there exists a map
$\Hat{P}(f)\colon\Hat{P}(X)\to\Hat{P}(Y)$. A systematic study of the
functor $\Hat{P}$ can be found in \cite{b1} and \cite{b2}. We also
consider the subspace $P_\beta(X)\subset\Hat{P}(X)$ consisting of
all measures $\mu$ such that $supp\hbox{~}\mu$ is compact.

This paper is devoted to some properties of Milyutin maps between
metrizable spaces. We say that $f\colon X\to Y$ is a {\em Milyutin
map} if there exists a map $g\colon Y\to \Hat{P}(X)$ such that
$supp\hbox{~}g(y)\subset f^{-1}(y)$ for every $y\in Y$. Such $g$ is
called a choice map associated with $f$. According to \cite[Theorem
3.15]{b2}, for any metrizable $X$ there exists a barycentric map
$b_{\Hat{P}(X)}\colon\Hat{P}(\Hat{P}(X))\to\Hat{P}(X)$ such that
$b_{\Hat{P}(X)}(\nu)=\nu$ for all $\nu\in\Hat{P}(X)$. Hence, if $g$
is a choice map associated with $f$, then the map
$b_{\Hat{P}(X)}\circ\Hat{P}(g)\colon\Hat{P}(Y)\to\Hat{P}(X)$ is a
right inverse of $\Hat{P}(f)$. Consequently, $f$ is a Milyutin map
if and only if $\Hat{P}(f)$ admits a right inverse.

Our first principal result concerns the question when $\Hat{P}(f)$
is a soft map. Recall that a map $f\colon X\to Y$ is soft if for any
space $Z$ and its closed subset $A$ and any maps $g\colon Z\to Y$,
$h\colon A\to X$ with $(f\circ h)|A=g$ there exists a map
$\bar{g}\colon Z\to X$ such that $\bar{g}$ extends $h$ and
$f\circ\bar{g}=g$. It is easily seen that every soft map is
surjective and open.

\begin{thm}
Let $f\colon X\to Y$ be a surjective open map between complete
spaces. Then $\Hat{P}(f)\colon\Hat{P}(X)\to\Hat{P}(Y)$ is a soft
map.
\end{thm}

The particular cases of Theorem 1.1 when both $X$ and $Y$ are either
compact or separable were established in \cite{f} and \cite{b3},
respectively.

Since any soft map admits a right inverse, a map $f$ satisfying the
hypotheses of Theorem 1.1 is a  Milyutin map. We apply Theorem 1.1
to obtain some results about atomless and exact Milyutin maps
introduced in \cite{rs}. If $f\colon X\to Y$ is a Milyutin map and
there exists a choice map $g$ such that $supp\hbox{~}
g(y)=f^{-1}(y)$ (resp., $g(y)$ is an atomless measure on $f^{-1}(y)$
for each $y\in Y$, i.e. $g(y)(\{x\})=0$ for all $x\in f^{-1}(y)$),
then $f$ is said to be an {\em exact} (resp., {\em atomless})
Milyutin map. It was established in \cite{rs} that, in the realm of
Polish spaces $X$ and $Y$, $f$ is exact Milyutin if and only if it
is open. The classes of atomless exact Milyutin maps and atomless
Milyutin maps between Polish spaces were characterized in
\cite[Theorem 1.6]{at}. The first class consists of all open maps
possessing perfect fibers (i.e., without isolated points)
\cite[Theorem 1.6]{at}, and the second one of all maps $f\colon X\to
Y$ such that for some  Polish space $X_0\subset X$ the restriction
$f_0=f|X_0\colon X_0\to Y$ is an open surjection whose fibers are
perfect \cite[Theorem 1.7]{at}.

Next theorem is a non-separable analogue of \cite[Theorem 1.7]{at}.

\begin{thm}
A continuous surjection $f\colon X\to Y$ of complete spaces is an
atomless Milyutin map if and only if there exists a complete
subspace $X_0\subset X$ such that $f_0=f|X_0\colon X_0\to Y$ is an
open surjection and all fibers of $f_0$ are perfect sets. Moreover,
for any such $f$ there exists a map $f^*\colon
P_\beta(Y)\to\Hat{P}(X)$ such that any $f^*(\mu)$ is atomless and
$\Hat{P}(f)\big(f^*(\mu)\big)=\mu$, $\mu\in P_\beta(Y)$.
\end{thm}

We do not know whether under the hypotheses in Theorem 1.2 there
exists a map $f^*\colon\Hat{P}(Y)\to\Hat{P}(X)$ such that each
$f^*(\mu)$ is atomless and $\Hat{P}(f)\big(f^*(\mu)\big)=\mu$,
$\mu\in\Hat{P}(Y)$. But if we are interested in atomless maps
defined on $Y$, we have the following:

\begin{thm}
Every open surjection $f\colon X\to Y$ with perfect fibers is a
densely atomless Milyutin map provided $X$ and $Y$ are complete
spaces.
\end{thm}

Here, a Milyutin map $f\colon X\to Y$ is {\em densely atomless}  if
$$\{g\in Ch_f(Y,X):g(y)\hbox{~}\mbox{is atomless for all}\hbox{~}y\in Y\}$$
is a dense $G_\delta$-set in the space $Ch_f(Y,X)$ of all choice
maps associated with $f$ equipped with the source limitation
topology. A few words about this topology. If $(X,d)$ is a bounded
(complete) metric space, then there exists a (complete) metric
$\Hat{d}$ on $\Hat{P}(X)$ generating its topology and extending $d$,
see \cite{b2}. Then $Ch_f(Y,X)$ is a subspace of the function space
$C(Y,\Hat{P}(X))$ with the source limitation topology whose local
base at a given $h\in C(Y,\Hat{P}(X))$ consists of all sets
$$B_{\Hat{d}}(h,\alpha)=\{g\in C(Y,\Hat{P}(X)):\Hat{d}(g(y),h(y))<\alpha(y)\hbox{~}\mbox{for all}\hbox{~}y\in Y\},$$ where
$\alpha$ is a continuous map from $Y$ into $(0,\infty)$. It is well
known that this topology does not depend on the metric $\Hat{d}$ and
it has the Baire property in case $\Hat{P}(X)$ is complete.
Similarly, $f$ is said to be {\em densely exact} provided the set
$$\{g\in Ch_f(Y,X):supp\hbox{~}g(y)=f^{-1}(y)\hbox{~}\mbox{for every}\hbox{~}y\in Y\}$$
is a dense and $G_\delta$-set in $Ch_f(Y,X)$. When $f$ is both
densely atomless and densely exact, it is called densely exact
atomless.

\begin{thm}
Let $f\colon X\to Y$ be an open surjection of complete spaces and
$\pi\colon X\to M$ a map into a separable space $M$. Then all choice
maps $h\in Ch_f(Y,X)$ such that $\pi(supp\hbox{~}h(y))$ is dense in
$\pi(f^{-1}(y))$ for every $y\in Y$ form a dense $G_\delta$-set in
$Ch_f(Y,X)$.
\end{thm}

It is interesting whether in Theorem 1.4 one can substitute the
phrase "$\pi(supp\hbox{~}h(y))$ is dense in $\pi(f^{-1}(y))$" by
"$\pi(supp\hbox{~}h(y)=\pi(f^{-1}(y))$".

Next corollary is a parametrization of the Parthasarathy \cite{pa}
result that perfect Polish spaces admit atomless measures. It also
provides a partial answer of the question \cite{at} whether any open
surjection $f$ of complete spaces is an exact atomless Milyutin map
provided all fibers of $f$ are perfect Polish spaces.

\begin{cor}
Let $f\colon X\to Y$ be an open and closed surjection of complete
 spaces such that all fibers of $f$ are separable $($and
perfect$)$. Then  $f$ is densely exact $($atomless$)$ Milyutin map.
\end{cor}

Finally, we generalize \cite[Corollary 1.4]{rs} and \cite[Corollary
1.9]{at} as follows (below a continuous set-valued map means a map
which is both lower and upper semi-continuous):

\begin{cor}
Let $X$ and $Y$ be complete  spaces and $\Phi\colon Y\to X$ a
continuous set-valued map such that all values $\Phi(y)$ are closed
separable subsets of $X$. Then there exists a map $h\colon Y\to
\Hat{P}(X)$ such that $supp\hbox{~}h(y)=\Phi(y)$ for every $y\in Y$.
If, in addition, all $\Phi(y)$ are perfect sets, the map $h$ can be
chosen so that every $h(y)$ is atomless.
\end{cor}

\noindent\textbf{Acknowledgements:} The author wishes to thank the
referee for his/her valuable remarks and suggestions which
significantly improved the paper.

\section{Preliminaries}

In this section we provide some preliminary results and establish
the proof of Theorem 1.1.

Probability Radon measures on a complete  space $X$ can be described
as positive linear functionals $\mu$ on $C^*(X)$ such that
$||\mu||=1$ and $\lim\mu(h_\alpha)=0$ for any decreasing net
$\{h_\alpha\}\subset C^*(X)$ which pointwisely converges to 0, see
\cite{vr}. Under this interpretation, $supp\hbox{~}\mu$ coincides
with the set of all $x\in X$ such that for every neighborhood $U_x$
of $x$ in $X$ there exists $\varphi\in C^*(X)$ such that
$\varphi(X\backslash U_x)=0$ and $\mu(\varphi)\neq 0$. This
representation of $supp\hbox{~}\mu$ easily implies that the
set-valued map $supp\colon \Hat{P}(X)\to X$ (assigning to each $\mu$
its support) is lower semi-continuous, i.e., $\{\mu\in
\Hat{P}(X):supp\hbox{~}\mu\cap U\neq\varnothing\}$ is open in
$\Hat{P}(X)$ for any open $U\subset X$. For every closed $F\subset
X$, we have $\mu (F)=\inf\{\mu(\varphi): \varphi\in C(F)\}$ (see for
example \cite{ff} in case $X$ is compact), where $C(F)=\{\varphi\in
C^*(X): 0\leq\varphi\leq 1\hbox{~}\mbox{and}\hbox{~}\varphi(F)=1\}$.

According to \cite{b3}, any compatible (complete) metric $d$ on $X$
generates a compatible (complete) metric $\Hat{d}$ on $\Hat{P}(X)$
such that
$$\Hat{d}\big(t\mu+(1-t)\mu', t\nu+(1-t)\nu'\big)\leq
t\Hat{d}(\mu,\mu')+(1-t)\Hat{d}(\nu,\nu')$$ for all $t\in [0,1]$ and
$\mu,\mu',\nu,\nu'\in \Hat{P}(X)$. It is easily seen that every ball
(open or closed) with respect to $\Hat{d}$ is convex.

Let $A_\varepsilon(X)$ denote the set of all $\mu\in \Hat{P}(X)$
such that $\mu(\{x\})\geq\varepsilon$ for some $x\in
supp\hbox{~}\mu$. For any closed $K\subset X$ there exists a closed
embedding $i\colon \Hat{P}(K)\to \Hat{P}(X)$ defined by
$i(\nu)(h)=\nu(h|K)$ for all $\nu\in \Hat{P}(K)$ and $h\in C^*(X)$.
Everywhere below we identify $\Hat{P}(K)$ with the set
$i\big(\Hat{P}(K)\big)=\{\mu\in \Hat{P}(X):supp\hbox{~}\mu\subset
K\}$ which is closed in $\Hat{P}(X)$.

\begin{lem}
Let $X$ be a complete space, $K$ a perfect closed subset of $X$ and
$G$ a convex open subset of $\Hat{P}(K)$. Then for every
$\varepsilon>0$ we have:
\begin{enumerate}
\item $A_\varepsilon(X)$ is a closed subset of
$\Hat{P}(X)$;
\item $A_\varepsilon(X)\cap\overline{G}$ is a nowhere dense set in the closure $\overline{G}$.
\end{enumerate}
\end{lem}

\begin{proof}
$(1)$ Since $\Hat{P}(X)$ is metrizable, it suffices to check that
$\mu_0=\lim\mu_n\in A_\varepsilon(X)$ for every convergent sequence
$\{\mu_n\}_{n\geq 1}$ in $\Hat{P}(X)$ with $\{\mu_n\}\subset
A_\varepsilon(X)$. To this end, let $H$ be the closure in $X$ of the
set $\bigcup_{n\geq 0}supp\hbox{~}\mu_n$. Because every $\mu\in
\Hat{P}(X)$ has a separable support, $H$ is a Polish subset of $X$.
Considering all $\mu_n$, $n\geq 0$, as elements of $\Hat{P}(H)$, we
have that the sequence $\{\mu_n\}_{n\geq 1}$ is contained in
$A_\varepsilon(H)$ and converges to $\mu_0$. Therefore, by
\cite[Theorem 8.1]{pa}, $\mu_0\in A_\varepsilon(H)$. Consequently,
there exists $x_0\in H$ with $\mu_0(\{x_0\})\geq\varepsilon$.
Therefore, $A_\varepsilon(X)$ is closed in $\Hat{P}(X)$.

$(2)$ Since $A_\varepsilon(K)=A_\varepsilon(X)\cap\Hat{P}(K)$, it
suffices to show that $A_\varepsilon(K)$ is nowhere dense in
$\Hat{P}(K)$. Suppose $A_\varepsilon(K)$ contains an open subset $W$
of $\Hat{P}(K)$ and let $P_\omega(K)$ be the set of all $\mu\in
\Hat{P}(K)$ having a finite support. Since $P_\omega(K)$ is dense in
$\Hat{P}(K)$, there exists
$\mu_0=\sum_{i=1}^{i=k}\lambda_i\delta_{x_i}\in P_\omega(K)\cap W$.
Here, $\delta_{x_i}$ denotes Dirac's measures at $x_i$ and
$\lambda_i=\mu_0(\{x_i\})$. Moreover, $\lambda_i\geq\varepsilon$ for
at least one $i$. For each $i\leq k$ and $n\geq 1$ choose a
neighborhood $V_i\subset K$ of $x_i$ and $n$ different points
$x_{i(1)},..,x_{i(n)}\in V_i$ such that the family $\{V_i:1\leq
i\leq k\}$ is disjoint and $dist(x_i,x_{i(j)})\leq 1/n$ for all
$1\leq j\leq n$. This can be done because $K$ is perfect, so every
neighborhood of $x_i$ contains infinitely many points. Consider now
the measures
$\displaystyle\mu_n=\sum_{i=1}^{i=k}\sum_{j=1}^{j=n}\frac{\lambda_i\delta_{x_{i(j)}}}{n}$.
Since $\lim\mu_n=\mu_0$, there exists $n_0$ such that $\mu_n\in W$
for all $n\geq n_0$. Consequently, for every $n\geq n_0$ there
exists $i\leq k$ with $\lambda_i/n\geq\varepsilon$, a contradiction.
\end{proof}

\begin{lem}\label{soft}
Let $f\colon X\to Y$ be an open surjection between complete spaces
such that $dim Y=0$. Then $\Hat{P}(f)\colon \Hat{P}(X)\to
\Hat{P}(Y)$ is a soft map.
\end{lem}

\begin{proof}
According to Theorem 1.3 from \cite{b3}, it suffices to show that
$f$ is everywhere locally invertible. The last notion is defined as
follows: for any space $Z$, a point $a\in Z$, a map $g\colon Z\to Y$
and an open set $U\subset X$ with $g(a)\in f(U)$ there exist a
neighborhood $V$ of $a$ in $Z$ and a map $h\colon V\to U$ such that
$f\circ h=g|V$. Obviously, $f$ is everywhere locally invertible
provided it satisfies the following condition:

\begin{enumerate}
\item[(*)]
For any open $U\subset X$ and $a\in f(U)$ there exists a map
$g\colon V\to U$ with $V$ being a neighborhood of $a$ in $Y$ such
that $f(g(y))=y$ for all $y\in V$.
\end{enumerate}

To show $f$ satisfies $(*)$, fix an open set $U\subset X$ and $a\in
f(U)$. Since $f$ is open, the set $V=f(U)\subset Y$ is also open and
the set-valued map $\Phi\colon V\to U$, $\Phi(y)=f^{-1}(y)\cap U$,
is lower semi-continuous with closed values. Moreover, $U$ admits a
complete metric because $X$ is complete. Then, by the
$0$-dimensional selection theorem of Michael \cite{m3}, $\Phi$ has a
continuous selection $g$. Obviously, $g$ is as required.
\end{proof}

\noindent\textbf{Proof of Theorem 1.1.} First, let us show that
$\hat{f}=\Hat{P}(f)|\Hat{P}(f)^{-1}(Y)$ is everywhere locally
invertible. It suffices to show that $\hat{f}$ satisfies condition
$(*)$ from Lemma~\ref{soft}. Suppose that $U\subset
\Hat{P}(f)^{-1}(Y)$ is open and $y_0\in \hat{f}(U)$. We need to find
a map $\alpha\colon V\to U$, where $V$ is a neighborhood of $y_0$ in
$Y$, such that $\hat{f}(\alpha(y))=y$ for every $y\in V$. To this
end, choose a $0$-dimensional complete space $Z$ and a perfect
Milyutin map $g\colon Z\to Y$, see \cite{ch} (recall that a map is
perfect if it is closed and has compact fibers). Next, consider the
pull-back $T=\{(z,x)\in Z\times X:g(z)=f(x)\}$ of $Z$ and $X$ with
respect to the maps $g$ and $f$, and let $p_f\colon T\to Z$,
$p_g\colon T\to X$ be the corresponding projections. Since $f$ is
open, so is $p_f$. For any $y\in Y$ we have
$p_f^{-1}\big(g^{-1}(y)\big)=p_g^{-1}\big(f^{-1}(y)\big)=g^{-1}(y)\times
f^{-1}(y)$. Since $g$ is Milyutin, there exists a map $g^*\colon
Y\to \Hat{P}(Z)$ such that $supp\hbox{~}g^*(y)\subset g^{-1}(y)$ for
all $y\in Y$. Let $\hat{p_f}=\Hat{P}(p_f)\colon \Hat{P}(T)\to
\Hat{P}(Z)$ and $\hat{p_g}=\Hat{P}(p_g)\colon \Hat{P}(T)\to
\Hat{P}(X)$. Take an open set $G\subset \Hat{P}(X)$ with $G\cap
\Hat{P}(f)^{-1}(Y)=U$ and let $W=\hat{p_g}^{-1}(G)$. Pick $\mu^*\in
G\cap\Hat{P}(f^{-1}(y_0)$ and let $\nu_0=\mu_0\times\mu^*$ be the
product measure, where $\mu_0=g^*(y_0)$. Obviously,
$\nu_0\in\Hat{P}(g^{-1}(y_0)\times f^{-1}(y_0))\subset\Hat{P}(T)$.
Moreover, $\hat{p_f}(\nu_0)=\mu_0$ and $\nu_0\in W$ because
$\hat{p_g}(\nu_0)=\mu_*\in G$.

Now we can complete the proof that $\hat{f}$ is everywhere locally
invertible. Let $g_0\colon\{y_0\}\to\Hat{P}(T)$ be the constant map
$g_0(y_0)=\nu_0$. Since $\hat{p_f}(\nu_0)=g^*(y_0)$ and, by Lemma
2.2, the map $\hat{p_f}$ is soft, there exists a map $\theta\colon
Y\to\Hat{P}(T)$ extending $g_0$ such that
$\hat{p_f}\circ\theta=g^*$. Obviously, $V=\theta^{-1}(W)$ is a
neighborhood of $y_0$, and define $\alpha=\hat{p_g}\circ\theta$.
Since for any $y\in V$ we have $\hat{p_f}(\theta(y))=g^*(y)$,
$p_f(supp\hbox{~}\theta(y))=supp\hbox{~}g^*(y)\subset g^{-1}(y)$ and
$supp\hbox{~}\theta(y)\subset g^{-1}(y)\times f^{-1}(y)$. So,
$supp\hbox{~}\alpha(y)=p_g(supp\hbox{~}\theta(y))\subset f^{-1}(y)$.
Consequently, $\hat{f}(\alpha(y))=y$. Moreover, $\alpha(y)\in U$ for
all $y\in V$.

Since $\hat{f}$ is everywhere locally invertible, by \cite[Theorem
1.3]{b3}, the map $\Hat{P}(\hat{f})\colon \Hat{P}(\Hat{Y})\to
\Hat{P}(Y)$ is soft, where $\Hat{Y}=\hat{f}^{-1}(Y)$. Moreover,
$\Hat{P}(X)\subset\Hat{P}(\Hat{Y})\subset\Hat{P}(\Hat{P}(X))$
because $X\subset \Hat{Y}\subset\Hat{P}(X)$. Therefore the following
diagram
$$
\begin{CD}
\Hat{P}(\Hat{Y})@>{b_{\Hat{P}}}>>\Hat{P}(X)\cr
@V{\Hat{P}(\hat{f})}VV @VV{\Hat{P}(f)}V\cr
\Hat{P}(Y)@>{i_{\Hat{P}(Y)}}>>\Hat{P}(Y)\cr
\end{CD}
$$

\noindent is commutative. Here, $b_{\Hat{P}}$ denotes the
restriction $b_{\Hat{P}(X)}|\Hat{P}(\Hat{Y})$ of the barycentric map
$b_{\Hat{P}(X)}\colon\Hat{P}(\Hat{P}(X))\to\Hat{P}(X)$, see
\cite{b2}, and $i_{\Hat{P}(Y)}$ is the identity on $\Hat{P}(Y)$.
Since $b_{\Hat{P}}$ retracts each $\Hat{P}(\hat{f})^{-1}(\mu)$ onto
$\Hat{P}(f)^{-1}(\mu)$, $\mu\in\Hat{P}(Y)$, and $\Hat{P}(\hat{f})$
is soft, we finally obtain that $\Hat{P}(f)$ is also soft. The proof
is completed.

\section{Atomless Milyutin maps}

In this section we provide the proofs of Theorems 1.2 and 1.3.

\noindent \textbf{Proof of Theorem 1.2}. Suppose that $f\colon X\to
Y$ is a surjective atomless Milyutin map with $X$ and $Y$ complete
spaces. Then there exists a choice map $h\colon Y\to \Hat{P}(X)$
associated with $f$ such that $h(y)$ is an atomless measure for all
$y\in Y$. Let $X_0=\bigcup\{supp\hbox{~}h(y):y\in Y\}$ and
$f_0=f|X_0$. Since $f_0^{-1}=supp\circ h$ is lower semi-continuous,
$f_0$ is open. Hence, by \cite[Theorem 3.6]{at}, $X_0$ is complete.
Moreover, all $f_0^{-1}(y)$ are perfect sets because $h(y)$  are
atomless measures.

For the other implication, assume that $f\colon X\to Y$ is a
surjection between complete  spaces and there exists a complete
subspace $X_0\subset X$ such that $f_0=f|X_0$ is an open surjection
possessing perfect fibers.  Considering $X_0$ and $f_0|X_0$, we may
suppose that $f$ is open and all of its fibers $f^{-1}(y)$, $y\in
Y$, are perfect sets. Then, by Theorem 1.1,  $f$ is Milyutin because
$\Hat{P}(f)$ has a right inverse as a soft map. To show $f$ is
atomless, as in the proof of Theorem 1.1 take a $0$-dimensional
complete space $Z$ and a perfect Milyutin map $g\colon Z\to Y$.
Since $g$ is Milyutin, there exists a map
$g^*\colon\Hat{P}(Y)\to\Hat{P}(Z)$ such that
$\Hat{P}(g)\big(g^*(\mu)\big)=\mu$ for all $\mu\in\Hat{P}(Y)$. By
Theorem 1.1, $\Hat{P}(f)$ is open (as a soft map). Hence,
$\hat{f}\colon\Hat{P}(f)^{-1}(Y)\to Y$ is also open (as a
restriction of an open map onto a preimage-set). So, the set-valued
map $\Phi\colon Z\to\Hat{P}(f)^{-1}(Y)$,
$\Phi(z)=\hat{f}^{-1}(g(z))$, is lower semi-continuous. Actually,
$\Phi(z)=\Hat{P}(f^{-1}(g(z)))$ for every $z\in Z$.  Let $A_n$,
$n\geq 1$, be the set of all $\mu\in\Hat{P}(X)$ such that
$\mu(\{x\})\geq 1/n$ for some point $x\in supp\hbox{~}\mu$. Since
the fibers $f^{-1}(y)$ are perfect sets, by Lemma 2.1, $A_n$ are
closed in $\Hat{P}(X)$ and all intersections
$A_n\cap\Hat{P}(f^{-1}(y))$ are nowhere dense in
$\Hat{P}(f^{-1}(y))$, $y\in Y$. Then, by \cite[Theorem 1.2]{gv},
$\Phi$ admits a selection $\theta\colon Z\to\Hat{P}(f)^{-1}(Y)$ such
that $\theta(z)\in\Phi(z)\backslash\bigcup_{n=1}^{\infty}A_n$, $z\in
Z$. This means that each measure $\theta(z)\in\Hat{P}(f^{-1}(g(z)))$
is atomless. The selection $\theta$ generates a regular operator
$u\colon C^*(X)\to C^*(Z)$, $u(\phi)(z)=\theta(z)(\phi)$ for all
$\phi\in C^*(X)$ and $z\in Z$. Finally, for every $\mu\in
P_\beta(Y)$ let $f^*(\mu)\in\Hat{P}(X)$ be the measure defined by
$f^*(\mu)(\phi)=g^*(\mu)(u(\phi))$, $\phi\in C^*(X)$. It is easily
seen that this definition is correct (i.e., $f^*(\mu)\in\Hat{P}(X)$)
and $f^*\colon P_\beta(Y)\to\Hat{P}(X)$ is a continuous map.

Let us show that $\Hat{P}(f)\big(f^*(\mu)\big)=\mu$ for every
$\mu\in P_\beta(Y)$. It suffices to prove that $f^*(\mu)(\alpha\circ
f)=\mu(\alpha)$ for any $\alpha\in C^*(Y)$. And this is really true
because $\phi=\alpha\circ f$ is the constant $\alpha(y)$ on each set
$f^{-1}(y)$, $y\in Y$. So, $u(\phi)(z)=\theta(z)(\phi)=\alpha(y)$
for any $z\in g^{-1}(y)$. Thus, $u(\phi)=\alpha\circ g$ and
$f^*(\mu)(\alpha\circ f)=g^*(\mu)(\alpha\circ g)$. Finally, since
$\Hat{P}(g)\big(g^*(\mu)\big)=\mu$, we have $g^*(\mu)(\alpha\circ
g)=\mu(\alpha)$.

So, it remains to prove only that every $f^*(\mu)$, $\mu\in
P_\beta(Y)$, is an atomless measure. To this end, fix $\mu_0\in
P_\beta(Y)$, $x_0\in supp\hbox{~}f^*(\mu_0)$ and $\eta>0$. It
suffices to find a function $\phi_0\in C^*(X)$ with $0\leq\phi_0\leq
1$ such that $\phi_0(x_0)=1$ and $f^*(\mu_0)(\phi_0)\leq\eta$. Since
$\theta(z)(\{x_0\})=0$, for every $z\in Z$ there exists $\phi_z\in
C^*(X)$ and a neighborhood $U_z$ of $z$ in $Z$ such that
$0\leq\phi_z\leq 1$, $\phi_z(x_0)=1$ and $\theta(z')(\phi_z)<\eta$
whenever $z'\in U_z$. Using the compactness of
$g^{-1}(supp\hbox{~}\mu_0)$ (recall that $\mu_0$ has a compact
support and $g$ is a perfect map), we find neighborhoods $U_{z(i)}$,
$i=1,..,k$, covering $g^{-1}(supp\hbox{~}\mu_0)$, and let
$\phi_0=\phi_{z(1)}\cdot\phi_{z(2)}\cdot..\cdot\phi_{z(k)}$. Then
$\phi_0$ is as required. Indeed, since
$\Hat{P}(g)\big(g^*(\mu_0)\big)=\mu_0$, $g^{-1}(supp\hbox{~}\mu_0)$
contains the support of $g^*(\mu_0)$. Consequently,
$g^*(\mu_0)(u(\phi_0))\leq\max\{u(\phi_0)(z):z\in
g^{-1}(supp\hbox{~}\mu_0)\}$. So, there exists $z_0\in
g^{-1}(supp\hbox{~}\mu_0)$ such that $g^*(\mu_0)(u(\phi_0))\leq
u(\phi_0)(z_0)$. Next, choose $j$ with $z_0\in U_{z(j)}$ and observe
that $\phi_0\leq\phi_j$ implies $u(\phi_0)(z_0)\leq
u(\phi_j)(z_0)=\theta(z_0)(\phi_j)$. Therefore,
$f^*(\mu_0)(\phi_0)\leq \theta(z_0)(\phi_j)<\eta$ because $z_0\in
U_{z(j)}$. The proof is completed.

\smallskip
\noindent\textbf{Proof of Theorem 1.3}. Take a $0$-dimensional
complete space $Z$, a perfect Milyutin map $g\colon Z\to Y$ and a
map $g^*\colon\Hat{P}(Y)\to\Hat{P}(Z)$ which is a right inverse of
$\Hat{P}(g)$. We equip $\Hat{P}(X)$ with a convex metric $\Hat{d}$,
and let $A_n$, $n\geq 1$, be the closed subsets of $\Hat{P}(X)$
considered in the proof of Theorem 1.2. We need to show that the set
$\mathcal{A}$ of all atomless choice maps form a dense
$G_\delta$-subset of $Ch_f(Y,X)$. Since each $A_n$ is closed in
$\Hat{P}(X)$, it is easily seen that the sets
$$\mathcal{U}_n=\{h\in Ch_f(Y,X):
h(y)\not\in A_n\hbox{~}\mbox{for all}\hbox{~}y\in Y\}$$ are open in
$Ch_f(Y,X)$ and $\mathcal{A}=\bigcap_{n\geq 1}\mathcal{U}_n$. To
prove that $\mathcal{A}$ is dense in $Ch_f(Y,X)$, fix $h\in
Ch_f(Y,X)$ and a function $\eta\colon Y\to (0,\infty)$. We are going
to find a map $h'\in\mathcal{A}$ such that
$\Hat{d}(h(y),h'(y))\leq\eta(y)$ for all $y\in Y$.

Denote by $B(h(g(z)),\eta(g(z)))$ the open ball in $\Hat{P}(X)$
(with respect to $\Hat{d}$) which is centered at $h(g(z))$ and has a
radius $\eta(g(z))$. Define the set-valued map $\Phi\colon
Z\to\Hat{P}(X)$, $\Phi(z)=\overline{\Hat{P}(f^{-1}(g(z)))\cap
B(h(g(z)),\eta(g(z)))}$.  This is a convex and closed-valued map
because any ball in $\Hat{P}(X)$ with respect to $\Hat{d}$ is
convex. Since $\hat{f}=\Hat{P}(f)|\big(\Hat{P}(f)^{-1}(Y)\big)$ is
open (as a soft map, see Theorem 1.1), the set-valued map $z\mapsto
\Hat{P}(f)^{-1}(g(z))$ is lower semi-continuous. Hence, by
\cite[Proposition 2.5]{m1}, so is $\Phi$. Moreover, each $\Phi(z)$
is the closure of the convex open set $\Hat{P}(f^{-1}(g(z)))\cap
B(h(g(z)),\eta(g(z)))$ in $\Hat{P}(f^{-1}(g(z)))$. Hence, according
to Lemma 2.1, $A_n\cap\Phi(z)$, $n\geq 1$, are nowhere dense sets in
$\Phi(z)$ for every $z\in Z$. Then, by \cite[Theorem 1.2]{gv},
$\Phi$ has  a continuous selection $\theta\colon Z\to\Hat{P}(X)$
avoiding the set $\bigcup_{n=1}^{\infty}A_n$, i.e., with
$\theta(z)\in\Phi(z)\backslash\bigcup_{n=1}^{\infty}A_n$ for every
$z\in Z$. Following the notations from the proof of Theorem 1.2, we
extend $\theta$ to a map $\bar{\theta}\colon
P_\beta(Z)\to\Hat{P}(X)$ by $\bar{\theta}(\nu)(\phi)=\nu(u(\phi))$,
$\phi\in C^*(X)$. Now let $h'\colon Y\to\Hat{P}(X)$ be the
composition $\bar{\theta}\circ g^*$. It follows from the proof of
Theorem 1.2 that $h'(y)$ is atomless and
$h'(y)\in\Hat{P}(f^{-1}(y))$ for all $y\in Y$. So,
$h'\in\mathcal{A}$.

It remains to show that $\Hat{d}(h(y),h'(y))\leq\eta(y)$, $y\in Y$.
To this end, we fix $y\in Y$ and take a sequence $\{\nu_n\}\subset
P_\beta(g^{-1}(y))$ converging to $g^*(y)$ such that each $\nu_n$
has a finite support. It is easily seen that if
$\nu=\sum_{i=1}^{i=k}t_i\delta_{z(i)}\in P_\beta(g^{-1}(y))$ is a
measure with a finite support, then
$\bar{\theta}(\nu)=\sum_{i=1}^{i=k}t_i\theta(z(i))$. Since
$\Hat{d}(\theta(z(i)),h(y))\leq\eta(y)$ for all $i$ and the metric
$\Hat{d}$ is convex, we have
$\Hat{d}(\bar{\theta}(\nu),h(y))\leq\eta(y)$. In particular,
$\Hat{d}(\bar{\theta}(\nu_n),h(y))\leq\eta(y)$ for every $n$. This
implies that $\Hat{d}(h'(y),h(y))\leq\eta(y)$ because $h'(y)$ is the
limit of the sequence $\{\bar{\theta}(\nu_n)\}$.

\section{Exact Milyutin maps}

In this section the proofs of Theorem 1.4 and Corollaries 1.5-1.6
are established.

\begin{lem}
Let $U\subset X$ be a non-empty open set in a space $X$. Then the
set $\Hat{U}=\{\nu\in\Hat{P}(X):supp\hbox{~}\nu\cap
U\neq\varnothing\}$ is open convex and dense in $\Hat{P}(X)$.
\end{lem}

\begin{proof}
Since the support map $\nu\rightarrow supp\hbox{~}\nu$ is a lower
semi-continuous map, $\Hat{U}\subset\Hat{P}(X)$ is open. To show it
is dense, suppose there exists an open set
$W=\{\nu\in\Hat{P}(X):|\nu(\phi_i)-\nu_0(\phi_i)|<\varepsilon, 1\leq
i\leq k\}$ in $\Hat{P}(X)$ with
$W\subset\Hat{P}(X)\backslash\Hat{U}$, where $\phi_i\in C^*(X)$ and
$\varepsilon>0$. We can suppose that $\nu_0$ has a finite support
(recall that the measures with a finite support form a dense set in
$\Hat{P}(X)$). Let $\nu_0=\sum_{j=1}^{j=m}\lambda_j\delta_{x(j)}$
such that $\lambda_j>0$ and $\sum_{j=1}^{j=m}\lambda_j=1$. Then
$supp\hbox{~}\nu_0=\{x(j):1\leq j\leq m\}\subset X\backslash U$.
Now, let
$\nu'=\lambda_0\delta_{x(0)}+(\lambda_1-\lambda_0)\delta_{x(1)}+\sum_{j=2}^{j=m}\lambda_j\delta_{x(j)}$,
where $x_0\in U$ and $0<\lambda_0<\lambda_1$ such that
$\lambda_0|\phi_i(x_0)-\phi_i(x_1)|<\epsilon$ for every
$i=1,2,..,k$. The choice of $\lambda_0$ yields that $\nu'\in W$.
Consequently, $\nu'\not\in\Hat{U}$ and $supp\hbox{~}\nu'\subset
X\backslash U$. This contradicts $x_0\in U\cap supp\hbox{~}\nu'$.

To show $\Hat{U}$ is convex, it suffices to prove that
$supp\hbox{~}\big(t\nu_1+(1-t)\nu_2\big)=supp\hbox{~}\nu_1\cup
supp\hbox{~}\nu_2$ for any $\nu_1,\nu_2\in\Hat{P}(X)$ and any $t\in
(0,1)$. Obviously, $supp\hbox{~}\nu_1\cup supp\hbox{~}\nu_2\supset
supp\hbox{~}\big(t\nu_1+(1-t)\nu_2\big)$. Assume $x\in
supp\hbox{~}\nu_1$. Then for every neighborhood $V_x$ of $x$ there
exists a function $\phi_x\in C^*(X)$ with $\phi_x(X\backslash
V_x)=0$ and $\nu_1(\phi_x)\neq 0$. Since
$\nu_1(\phi_x)=\nu_1(\phi_x^+)-\nu_1(\phi_x^-)$, where $\phi_x^+$
and $\phi_x^-$ are the positive and negative parts of $\phi_x$, we
can suppose $\phi_x$ is non-negative. Then,
$\nu(\phi_x)\geq\nu_1(\phi_x)>0$ with $\nu=\nu=t\nu_1+(1-t)\nu_2$.
Hence, $x\in supp\hbox{~}\nu$ which completes the proof.
\end{proof}

\smallskip\noindent
\textbf{Proof of Theorem 1.4.} Choose a countable base $\{V_n:n\geq
1\}$ for the topology of $M$, and let
$B_n=\{\nu\in\Hat{P}(X):supp\hbox{~}\nu\cap
\pi^{-1}(V_n)=\varnothing\}$. By Lemma 4.1, each $B_n$ is closed in
$\Hat{P}(X)$. Let  $\mathcal{B}$ be the set of all maps $h\in
Ch_f(Y,X)$ such that $\pi(supp\hbox{~}h(y))$ is dense in
$\pi(f^{-1}(y))$ for any $y\in Y$. Obviously,
$\mathcal{B}=\bigcap_{n\geq 1}\mathcal{G}_n$, where
$\mathcal{G}_n=\{h\in Ch_f(Y,X): h(y)\not\in B_n\hbox{~}\mbox{for
all}\hbox{~}y\in Y\}$. It suffices to show that each $\mathcal{G}_n$
is open and dense in $Ch_f(Y,X)$ with respect to the source
limitation topology.

\textit{Claim $1$. Each  $\mathcal{G}_n$ is open in $Ch_f(Y,X)$.}

We can suppose that each $V_n$ is of the form
$V_n=g_n^{-1}(0,\infty)$ for some non-negative function $g_n\in
C^*(M)$. Then  $\nu\in B_n$ if and only if $\nu(g_n\circ\pi)=0$,
$n\geq 1$. Obviously the equality
$D_n(\mu,\mu')=\Hat{d}(\mu,\mu')+|\mu(g_n\circ\pi)-\mu'(g_n\circ\pi)|$,
where $\mu, \mu'\in\Hat{P}(X)$ and $\Hat{d}$ is a compatible metric
on $\Hat{P}(X)$, defines a compatible metric on $\Hat{P}(X)$ for
every $n\geq 1$. Given $h\in\mathcal{G}_n$ we consider the
continuous function $\alpha\colon Y\to (0,\infty)$,
$\alpha(y)=h(y)(g_n\circ\pi)/2$. We have
$B_{D_n}(h,\alpha)\subset\mathcal{G}_n$. Indeed, if $h'\in
B_{D_n}(h,\alpha)$, then $|h'(y)(g_n\circ\pi)-h(y)(g_n\circ\pi)|\leq
D_n(h(y),h'(y))<\alpha(y)$ for all $y\in Y$. The last inequality
implies $h'(y)(g_n\circ\pi)>\alpha(y)>0$, $y\in Y$. Hence,
$h'(y)\not\in B_n$ for all $y\in Y$. So, $h'\in\mathcal{G}_n$ which
completes the proof of Claim 1.

To show that any $\mathcal{G}_n$ is dense in $Ch_f(Y,X)$, we fix
$m\geq 1$, $h\in Ch_f(Y,X)$ and a function $\eta\colon Y\to
(0,\infty)$. We are going to find a map $h'\in\mathcal{G}_m$ with
$\Hat{d}(h'(y),h(y))\leq\eta(y)$ for all $y\in Y$. To this end,
following the proof of Theorems 1.2 and 1.3, take a complete
0-dimensional space $Z$ and a perfect Milyutin map $g\colon Z\to Y$
with a right inverse $g^*\colon Y\to P_\beta(Z)$. We also consider
the lower semi-continuous convex and closed-valued map $\Phi\colon
Z\to\Hat{P}(X)$, $\Phi(z)=\overline{\Hat{P}(f^{-1}(g(z)))\cap
B(h(g(z)),\eta(g(z)))}$. According to Lemma 4.1,
$B_m\cap\Hat{P}(f^{-1}(g(z)))$ is a closed nowhere dense subsets of
$\Hat{P}(f^{-1}(g(z)))$ for every $z\in Z$. Hence, all
$B_m\cap\Phi(z)$ are closed and nowhere dense in $\Phi(z)$. Then, by
\cite[Theorem 1.2]{gv}, $\Phi$ has  a continuous selection
$\theta\colon Z\to\Hat{P}(X)$ such that
$\theta(z)\in\Phi(z)\backslash B_m$, $z\in Z$. As in the proof of
Theorem 1.3, let $h'\colon Y\to\Hat{P}(X)$ be the composition
$\bar{\theta}\circ g^*$, where $\bar{\theta}\colon
P_\beta(Z)\to\Hat{P}(X)$ is an extension of $\theta$ defined by
$\bar{\theta}(\nu)(\phi)=\nu(u(\phi))$, $\phi\in C^*(X)$. Following
the arguments from Theorem 1.3, we can show that
$\Hat{d}(h'(y),h(y))\leq\eta(y)$ for all $y\in Y$. Next claim
completes the proof of Theorem 1.4.

\textit{Claim $2$. $h'(y)\not\in B_m$ for any $y\in Y$.}

The proof of this claim is reduced to find a function $\phi_y\in
C^*(X)$ such that $\phi_y\big(X\backslash\pi^{-1}(V_m)\big)=0$ and
$h(y)(\phi_y)\neq 0$. Indeed, in such a case
$supp\hbox{~}h(y)\cap\pi^{-1}(V_m)\neq\varnothing$. Since
$\theta(z)\not\in B_m$ for all $z\in g^{-1}(y)$,
$supp\hbox{~}\theta(z)\cap\pi^{-1}(V_m)\neq\varnothing$.
 Consequently, for any $z\in
g^{-1}(y)$ there exists a function $\phi_z\in C^*(X)$ with
$\phi_z\big(X\backslash\pi^{-1}(V_m)\big)=0$ and
$\theta(z)(\phi_z)\neq 0$. Considering the positive or negative
parts of $\phi_z$, we may assume each $\phi_z\geq 0$. Next, use the
continuity of $\theta$ and the compactness of $g^{-1}(y)$ to find
finitely many points $z(i)\in g^{-1}(y)$, $i=1,2,..,k$, and
neighborhoods $U_{z(i)}$ such that $\theta(z)(\phi_{z(i)})>0$
provided $z\in U_{z(i)}$. Finally, let
$\phi_y=\sum_{i=1}^{i=k}\phi_{z(i)}$. Then
$\phi_y\big(X\backslash\pi^{-1}(V_m)\big)=0$ and
$u(\phi_y)(z)=\theta(z)(\phi_y)>0$ for any $z\in g^{-1}(y)$. So,
$h(y)(\phi_y)\geq\min\{u(\phi_y)(z):z\in g^{-1}(y)\}>0$ because
$g^{-1}(y)$ is compact. This completes the proof of the claim.

\smallskip
\noindent \textbf{Proof of Corollary 1.5.} Since $f$ is closed with
separable fibers, there exists a map $\pi\colon X\to Q$ such that
all restrictions $\pi|f^{-1}(y)$, $y\in Y$, are embeddings, see
\cite{bp}. Here, $Q$ is the Hilbert cube. Then, by Theorem 1.4 (with
$M$ replaced by $Q$), $f$ is densely exact. If, in addition, the
fibers of $f$ are perfect, both Theorems 1.3 and 1.4 imply that $f$
is densely exact atomless.

\smallskip \noindent
\textbf{Proof of Corollary 1.6.} Consider the graph
$G(\Phi)=\cup\{\{y\}\times\Phi(y):y\in Y\}\subset Y\times X$ of
$\Phi$ and the projection $f\colon G(\Phi)\to Y$. Since $\Phi$ is
continuous, $G(\Phi)$ is closed in $Y\times X$ and $f$ is both open
and closed. Then $G(\Phi)$ is a complete space.  Now, by Corollary
1.5, there exists a map $h'\colon Y\to\Hat{P}(G(\Phi))$ with each
$h'(y)\in\Hat{P}(f^{-1}(y)$ being exact measure. Therefore,
$supp\hbox{~}h'(y)=f^{-1}(y)$. Let $h=\Hat{P}(\pi)\circ h'$, where
$\pi\colon G(\Phi)\to X$ is the projection into $X$. Since $\pi$
embeds each $f^{-1}(y)$ onto $\Phi(y)$, $h$ is a map from $Y$ into
$\Hat{P}(X)$ such that $supp\hbox{~}h(y)=\Phi(y)$  for every $y\in
Y$. If $\Phi(y)$ are perfect sets, so are the fibers $f^{-1}(y)$,
and $h'$ can be chosen to be atomless and exact. In such a case $h$
is also atomless.

\smallskip \noindent
\textbf{Note added in proof.} Recently T. Banakh informed the author
that V. Bogachev and A. Kolesnikov \cite{bk} proved the following
result: The map $\Hat{P}(f)$ from Theorem 1.1 is open. This, in
combination with Michael's convex-valued selection theorem
\cite{m1}, provides another proof of Theorem 1.1.

\end{document}